\documentclass[11pt]{amsart}
\usepackage[latin1]{inputenc}
\usepackage{amsmath,amsthm,amsfonts,amssymb,amscd}
\usepackage{verbatim}
\usepackage{mathrsfs}
\usepackage{multicol}
\usepackage{tabularx}
\usepackage[usenames,dvipsnames]{color}
\usepackage{enumerate}
\usepackage{graphicx}
\usepackage{color,xcolor}
\usepackage{bbm}
\usepackage{young}
\usepackage[vcentermath]{youngtab}
\usepackage[all]{xy}
\usepackage{hyperref}

\theoremstyle{definition}
\newtheorem{thm}{Theorem}[subsection]
\newtheorem{prop}[thm]{Proposition}
\newtheorem{cor}[thm]{Corollary}

\newtheorem{lem}[thm]{Lemma}

\newtheorem{rem}[thm]{Remark}

\def\bro{\text{\Large$\rho$}}

\def\ord#1^#2{#1$^{\text{#2}}$}
\def\lie#1{\mathfrak{#1}}

\def\hlie#1{\hat{\mathfrak{#1}}}

\def\uqr#1^#2{\text{$U_q^{#2}(\lie #1)$}}

\def\uqhr#1^#2{\text{$U_q^{#2}(\hlie #1)$}}
\def\us#1^#2{\text{$U_{\xi}^{#2}(\lie #1)$}}
\def\ush#1^#2{\text{$U_{\xi}^{#2}(\hlie #1)$}}
\def\dus#1^#2{\text{$\dot{U}_{\xi}^{#2}(\lie #1)$}}
\def\dush#1^#2{\text{$\dot{U}_{\xi}^{#2}(\hlie #1)$}}

\def\wt{{\rm wt}}

\def\ch{{\rm ch}}

\def\opl_#1^#2{\text{\scriptsize$\bigoplus\limits_{\text{\normalsize$#1$}}^{\text{\normalsize$#2$}}$}}
\def\otm_#1^#2{\text{\scriptsize$\bigotimes\limits_{\text{\footnotesize$#1$}}^{\text{\footnotesize$#2$}}$}}

\def\qbinom#1#2{\begin{bmatrix} #1\\#2\end{bmatrix}}
\def\tqbinom#1#2{\text{$\left[\begin{smallmatrix} #1\\#2\end{smallmatrix}\right]$}}

\def\endd{\hfill$\diamond$}

\begin{document}
\title[Limits in Excellent Filtrations and Tensor Products]{Limits of Multiplicities in Excellent Filtrations\\ and Tensor Product Decompositions for\\ Affine Kac-Moody Algebras}
\author[Dijana Jakeli\'c and Adriano Moura]{Dijana Jakeli\'c and Adriano Moura}
\thanks{}

\address{Department of Mathematics and Statistics, University of North Carolina Wilmington, \vspace{-12pt}}
\address{601 S. College Road, Wilmington, NC, 28401-5970}
\email{jakelicd@uncw.edu}

\address{Department of Mathematics, University of  Campinas, Campinas - SP - Brazil, 13083-859.}
\email{aamoura@ime.unicamp.br}
\thanks{The work of the second author was partially supported by CNPq grant 304477/2014-1 and Fapesp grant 2014/09310-5.}

\begin{abstract}
We express the multiplicities of the irreducible summands of certain tensor products of irreducible integrable modules for an affine Kac-Moody algebra over a simply laced Lie algebra as sums of multiplicities in appropriate excellent filtrations (Demazure flags). As an application, we obtain expressions for the outer multiplicities of tensor products of two fundamental modules for $\widehat{\lie{sl}}_2$ in terms of partitions with bounded parts, which subsequently lead to certain partition identities.
\end{abstract}

\maketitle

\section{Introduction}
\numberwithin{equation}{section}

The tensor product of two integrable highest-weight modules for any symmetrizable Kac-Moody algebra decomposes into a  direct sum of integrable highest-weight modules with each summand having finite multiplicity, often called the outer multiplicity of that summand in the given tensor product. In principle, the outer multiplicities can be computed algorithmically \cite{kiwe,lit:LR,oss}. However, it is well-known that more direct descriptions lead to deep connections with combinatorics, number theory, and mathematical physics. For instance, in the case of affine Lie algebras this leads to proofs of Rogers-Ramanujan-type identities as well as partition identities (see \cite{fein,lepwil,mw:tp,mw:tp2} and references therein).

Unless the given Kac-Moody algebra is of finite type, all irreducible integrable modules (except the trivial) are infinite-dimensional. The Demazure modules provide a way of studying the irreducible integrable modules via certain finite-dimensional pieces. Given such an irreducible highest-weight module $V$ and an extremal weight vector $v\in V$, the associated Demazure module $D$  is the submodule for the standard positive Borel subalgebra generated by $v$. If $\hlie g$ is the non-twisted affine Kac-Moody algebra over a simple Lie algebra $\lie g$, one may ask whether $D$ is also a $\lie g$-submodule of $V$. This is the case exactly when the restriction of the weight of $v$ to the Cartan subalgebra of $\lie g$ is anti-dominant. In that case, $D$ is often referred to as a $\lie g$-stable Demazure module and it is also a graded submodule for the current algebra $\lie g[t]\subseteq\hlie g$. If $V$ is a module of level $\ell$, then $D$ is said to be a Demazure module of level $\ell$ (the level is the scalar by which the central element of $\hlie g$ acts on $V$). 

The finite-dimensional representation theory of $\lie g[t]$ is presently a very active research area motivated by, among other reasons, its application to the problem of understanding the structure of finite-dimensional representations of quantum affine algebras. It follows from the previous paragraph that $\lie g$-stable Demazure modules are examples of such representations of $\lie g[t]$. It is usual to denote by $D(\ell,\lambda,r)$ the $\lie g$-stable Demazure module of level $\ell$ having $w_0\lambda$ as the anti-dominant weight mentioned previously and the generating vector $v$ belonging to the $r$-th graded piece of $V$ (here $\lambda$ is a dominant weight for $\lie g$ and $w_0$ is the longest element of its Weyl group). When $r=0$, we simply write $D(\ell,\lambda)$ (the modules $D(\ell,\lambda,r)$ and $D(\ell,\lambda,r')$ are isomorphic as non-graded $\lie g[t]$-modules). The Demazure modules $D(\ell,\lambda)$ are quotients of the local Weyl module $W(\lambda)$, a certain universal module in the category of graded finite-dimensional $\lie g[t]$-modules. It turns out that we have an isomorphism of graded $\lie g[t]$-modules
\begin{equation*}
W(\lambda)\cong D(1,\lambda)
\end{equation*} 
when $\lie g$ is simply laced \cite{chlo:wfd,foli:weyldem}. Furtheremore, it was shown in \cite{naoi:weyldem} that, for any $\lie g$, $W(\lambda)$ admits a level-$1$ Demazure flag, i.e., a filtration whose subsequent quotients are isomorphic to Demazure modules of level $1$. This fact finally settled the problem of describing the character of $W(\lambda)$. Demazure flags are also known as {\it excellent filtrations} \cite{joseph03,joseph06}.

To explain how  Demazure flags are related to the outer multiplicities problem discussed in the first paragraph, it is useful to observe that, since the central element of $\hlie g$ is not in $\lie g[t]$, there are isomorphisms of graded $\lie g[t]$-modules between Demazure modules of different levels. In fact, if $\ell$ is sufficiently large (with respect to $\lambda$), then $D(\ell,\lambda)$ is isomorphic to the simple $\lie g$-module of highest weight $\lambda$ equipped with the trivial action of the positive graded pieces of $\lie g[t]$. Hence, it makes sense to ask, for any $\ell$, whether $D(\ell,\lambda)$ admits a level-$\ell'$ Demazure flag for $\ell'>\ell$.
The proof that Weyl modules admit level-$1$ Demazure flags given in \cite{naoi:weyldem} is based on a result providing an affirmative answer to this question in the case that $\lie g$ is simply laced. The proof of this result establishes the connection with the problem of computing outer multiplicities discussed earlier. 
Namely, the key ingredients in the proof were results from \cite{joseph03,joseph06} asserting that the tensor product of a singleton Demazure crystal $B_1$ with a general Demazure crystal $B_2$ can be decomposed as a disjoint union of Demazure crystals and that, under certain conditions, this decomposition is ``well-behaved'' with respect to global bases. The main idea behind the proof of the existence of a level-$(\ell+1)$ Demazure flag for $D(\ell,\lambda)$ in \cite{naoi:weyldem} was to use these results with $B_2$ being the crystal of $D(\ell,\lambda)$ and $B_1$ being the crystal associated to the one-dimensional Demazure module inside the basic representation $V(\Lambda_0)$ of $\hlie g$. Thus, the Demazure flag is essentially obtained from partial information about the decomposition of the tensor product
\begin{equation*}
V(\Lambda_0)\otimes V
\end{equation*}
where $V$ is the irreducible integrable $\hlie g$-module containing $D(\ell,\lambda)$. 

It has recently been shown that it may be possible to compute the multiplicities in Demazure flags of $\lie g$-stable Demazure modules through different methods, unrelated to the problem of outer multiplicities. This was proved for $\lie g=\lie{sl}_2$ in \cite{bcsv:dfcpmt,cssw:demflag}. Thus, the main idea of the present paper is to reverse the process of \cite{naoi:weyldem} and obtain the full information about the decomposition of $V(\Lambda_0)\otimes V$ from the knowledge of the multiplicities in the level-$(\ell+1)$ Demazure flags of all $\lie g$-stable Demazure modules contained in $V$.

We now describe our main result. All $\lie g$-stable Demazure modules inside $V$ are of the form $D(\ell,\mu,s)$. Let $$D(V)=\{(\mu,s): D(\ell,\mu,s)\subseteq V\}.$$
Note that, if $\ell>0$, $D(V)$ is an infinite set.
Suppose we want to compute the outer multiplicity of an irreducible module $W$ in $V(\Lambda_0)\otimes V$. The Demazure submodules of $W$ are of the form $D(\ell+1,\lambda,r)$ with $(\lambda,r)\in D(W)$. 
Denote by
\begin{equation*}
[D(\ell,\mu,s):D(\ell+1,\lambda,r)]
\end{equation*}
the multiplicity of $D(\ell+1,\lambda,r)$ in a Demazure flag for $D(\ell,\mu,s)$. 
Our main result (Theorem \ref{t:multrel}) states that the outer multiplicity of $W$ in $V(\Lambda_0)\otimes V$ is given by
\begin{equation*}
\sum_{(\lambda,r)\in D(W)}\ \max_{(\mu,s)\in D(V)}\    [D(\ell,\mu,s):D(\ell+1,\lambda,r)].
\end{equation*}
In particular, all summands are nonnegative and $\max\limits_{(\mu,s)\in D(V)}\    [D(\ell,\mu,s):D(\ell+1,\lambda,r)]\ne 0$ for finitely many $(\lambda,r)\in D(W)$. The above formula can then be generalized for tensor products of the form $V(\Lambda_i)\otimes V$ with $i$ in the orbit of $0$ under the action of the group of Dynkin diagram automorphisms of $\hlie g$ (see \eqref{e:multrelnot0}). 

In the case that $\lie g=\lie{sl}_2$ we describe the $\max$ above as the limit of the expressions for $[D(\ell,\mu,s):D(\ell+1,\lambda,r)]$ obtained in \cite{bcsv:dfcpmt,cssw:demflag}  (see Proposition \ref{c:sl2l}). Moreover, for $\ell=1$, we are able to express this limit as the cardinality of the set of bounded partitions of a certain nonnegative integer. Both this integer and the bound are given as functions of  $(\lambda,r)$  (see Proposition \ref{t:sl21}).  In particular, the subset of $D(W)$ which contributes to the summation is explicitly characterized. 
It is interesting to observe that, also for $\ell=1$, the outer multiplicities were recently calculated in \cite{mw:tp,mw:tp2} as the cardinality of a different set of partitions (for $\lie g$ of type $A$). By comparing the expressions from \cite{mw:tp,mw:tp2} with  Proposition \ref{t:sl21}, we obtain the partition identities \eqref{e:partrel}.  We believe that the results of \cite{mw:tp2} will also be helpful for carrying out the algorithm of \cite{joseph03} to obtain, in the spirit of \cite{naoi:weyldem},  expressions for the multiplicities of the level-$2$ Demazure flags of Weyl modules for type $A$ in terms of partitions associated to colored Young diagrams. We remark that, even in the $\lie{sl}_2$-case, this may lead to different expressions than those given in \cite{cssw:demflag}. We will address this in a future publication. We also point out that, for general $\ell$ and $\lie g=\lie{sl}_2$, the outer multiplicities were described in \cite{fein} in terms of alternating sums of partitions.

The paper is organized as follows. In Section \ref{s:main} we review the prerequisites about Weyl and Demazure modules that are necessary to state the main results: Theorem \ref{t:multrel} and Propositions \ref{c:sl2l} and \ref{t:sl21}. Theorem \ref{t:multrel} is proved in Section \ref{ss:multrel} following a review of the necessary facts on crystals, in particular the result of \cite{joseph03}  that we utilize.
In Section \ref{s:sl2} we specialize Theorem \ref{t:multrel} to the $\lie {sl}_2$-case and prove Propositions \ref{c:sl2l} and \ref{t:sl21}.

\section{The Main Results}\label{s:main}
\numberwithin{equation}{subsection}

\subsection{Weyl and Demazure Modules}
Given a Lie algebra $\lie a$ over $\mathbb C$, denote by $\lie a[t]$ the associated current algebra, i.e., $\lie a[t]=\lie a\otimes \mathbb C[t]$ with bracket $[x\otimes f,y\otimes g] = [x,y]\otimes (fg)$ for $x,y\in\lie a$ and $f,g\in\mathbb C[t]$. Set also  $\lie a[t]_+= \lie a\otimes t\mathbb C[t]$.

Let $\lie g$ be a finite-dimensional simple Lie algebra over $\mathbb C$, $R^+$  the set of positive roots with respect to a fixed triangular decomposition $\lie g=\lie n^-\oplus\lie h\oplus\lie n^+$, and $\Delta\subseteq R^+$  the corresponding set of simple roots. Fix a Chevalley basis $\{x_\alpha^\pm, h_\beta:\alpha\in R^+,\beta\in\Delta\}$ and set $h_\alpha=[x_\alpha^+,x_\alpha^-],\alpha\in R^+$. Let $P$ denote the weight lattice of $\lie g$ and $P^+$ the submonoid of dominant integral weights. We let $I$ denote the set of vertices of the Dynkin diagram of $\lie g$ and let $\alpha_i$ and $\omega_i, i\in I$, denote the corresponding simple root and fundamental weight, respectively. We often simplify notation and write $h_i=h_{\alpha_i}, x_i^\pm=x_{\alpha_i}^\pm, i\in I$.

Given $\lambda\in P^+$, the graded local Weyl module of highest weight $\lambda$ is the $\lie g[t]$-module $W(\lambda)$ defined by a generator $v$ satisfying the following defining relations:
\begin{equation}\label{e:W(lambda)}
\lie n^+[t]v = 0 = \lie h[t]_+v, \quad hv=\lambda(h)v, \quad (x_\alpha^-)^{\lambda(h_\alpha)+1}v = 0
\end{equation}
for all $h\in\lie h$ and all $\alpha\in\Delta$.  It is well-known that $W(\lambda)$ is finite-dimensional and any other finite-dimensional graded $\lie g[t]$-module generated by a vector satisfying the first three listed relations in \eqref{e:W(lambda)} is a quotient of $W(\lambda)$. Given $\ell\in\mathbb Z_{\ge 0}$, the $\lie g$-stable level-$\ell$ Demazure module $D(\ell,\lambda)$ is the quotient of $W(\lambda)$ by the submodule generated by
\begin{equation}\label{e:D(ell,lambda)}
\{(x_\alpha^-\otimes t^{p_\alpha}) v: \alpha\in R^+\}\cup \{(x_\alpha^-\otimes t^{p_\alpha-1})^{m_\alpha+1} v: \alpha\in R^+ \text{ s.t. } m_\alpha<\ell r^\vee_\alpha\},
\end{equation}
where $r^\vee_\alpha, p_\alpha$, and $m_\alpha$ are the integers  defined by
\begin{equation*}
r^\vee_\alpha = \begin{cases}
1,& \text{ if } \alpha \text{ is long,}\\
r^\vee,& \text{ if } \alpha \text{ is short},
\end{cases}
\quad\text{and}\quad \lambda(h_\alpha) = (p_\alpha-1)\ell r^\vee_\alpha + m_\alpha, \quad 0<m_\alpha\le \ell r^\vee_\alpha,
\end{equation*}
where $r^\vee$ is the lacing number of $\lie g$. The reason for the terminology ``$\lie g$-stable Demazure module'' will be explained by Lemma \ref{l:gstdem} and \eqref{e:cvdemrel}. 
It was proved in \cite{chlo:wfd} for $\lie g$ of type A and in \cite{foli:weyldem} for simply laced $\lie g$ that
\begin{equation}\label{e:sliso}
W(\lambda)\cong D(1,\lambda) \quad\text{for all}\quad \lambda\in P^+.
\end{equation}
For general $\lie g$ this is false for most values of $\lambda$ (see however, the paragraph preceding Theorem \ref{t:naoislf} below).

For $m\in\mathbb Z$, let $\tau_m$ be the functor that shifts degree by $m$. More precisely, if $V$ is a graded $\lie g[t]$-module and $V[k]$ is its $k$-th graded piece, then
\begin{equation*}
(\tau_m V)[k] = V[k-m].
\end{equation*}
Set
\begin{equation*}
D(\ell,\lambda,m) = \tau_m D(\ell,\lambda).
\end{equation*}
It is well known that 
\begin{equation}\label{e:isodeml}
 D(\ell,\lambda,m) \cong D(\ell,\mu,n) \qquad\Leftrightarrow\qquad (\lambda,m) = (\mu,n). 
\end{equation}

\subsection{Affine Algebras} The affine Kac-Moody algebra $\hlie g$ associated to $\lie g$ is the vector space $\hlie g=\lie g\otimes\mathbb C[t,t^{-1}]\oplus \mathbb Cc\oplus \mathbb Cd$ with bracket given by
\begin{equation*}
[x \otimes t^r, y \otimes t^s] = [x, y] \otimes t^{r+s} + r\ \delta_{r,-s}\ (x, y)\ c, \quad [c,\hlie g]=\{0\}, \quad\text{and}\quad [d,x\otimes t^r]=r\ x\otimes t^r
\end{equation*}
for any $x,y \in \lie g$ and $r,s \in \mathbb Z$, where $(\ ,\ )$ is the Killing form on $\lie g$. We identify $\lie g$ and $\lie g[t]$ with the obvious subalgebras of $\hlie g$. Set
\begin{equation*}
\hlie h = \lie h \oplus \mathbb C c \oplus \mathbb C d, \quad \hlie n^{+} = \lie n^{+} \oplus \lie g [t]_{+}, \quad\text{and}\quad \hlie b = \hlie n^{+} \oplus \hlie h.
\end{equation*}
Note that $\lie g[t]\oplus\mathbb Cc\oplus\mathbb Cd$ is a parabolic subalgebra of $\hlie g$ containing $\hlie b$. Identify $\lie h^*$ with the subspace $\{\lambda\in\hlie h^*:\lambda(c)=\lambda(d)=0\}$. Let $\Lambda_0,\delta\in\hlie h^*$ be defined by
\begin{equation*}
\Lambda_0(d)=0=\Lambda_0(\lie h), \quad \Lambda_0(c)=1,\quad \delta(c)=0=\delta(\lie h), \quad \delta(d)=1.
\end{equation*}
Also, set $\hat I=I\sqcup\{0\}, \alpha_0=\delta-\theta, h_0=c-h_\theta$, where $\theta$ is the highest root of $\lie g$, and, for $i\in I$, set
\begin{equation*}
\Lambda_i = \omega_i +\omega_i(h_\theta)\Lambda_0.
\end{equation*}
Then, $\Lambda_i(h_j)=\delta_{i,j}$ for all $i,j\in\hat I, \{h_i:i\in\hat I\}\cup\{d\}$ is a basis of $\hlie h$, $\hat\Delta=\{\alpha_i:i\in\hat I\}$ is the set of simple roots for $\hlie g$, and $\hat R^+ = R^+\cup\{\alpha+r\delta: \alpha\in R\cup\{0\}, r\in\mathbb Z_{>0}\}$ is the set of positive roots.

Set $\hat P = \{\lambda\in\hlie h^*: \lambda(h_i)\in\mathbb Z \text{ for all } i\in\hat I\}= \mathbb C\delta \oplus P\oplus \mathbb Z\Lambda_0$ and $\hat P^+ = \{\lambda\in\hat P: \lambda(h_i)\ge 0 \text{ for all } i\in\hat I\}$.
Given $\lambda\in\hat P$, the number $\lambda(c)$ is called the level of $\lambda$. Let
$\widehat{\mathcal W}$ denote the  affine Weyl group, which is generated by the simple reflections $s_i,i\in\hat I$. More precisely, $s_i$ is the linear endomorphism of $\hlie h^*$ defined by
\begin{equation*}
s_i\lambda = \lambda-\lambda(h_i)\alpha_i \quad\text{for all}\quad \lambda\in\hat P.
\end{equation*}
We denote by $\ell(w)$ the length of $w\in\widehat{\mathcal W}$, i.e., the length of the shortest expressions $w=s_{i_1}\cdots s_{i_\ell}, i_j\in\hat I, 1\le j\le\ell$. The Weyl group of $\lie g$ is the subgroup $\mathcal W$ of $\widehat{\mathcal W}$ generated by $s_i,i\in I$. Since $\alpha(c)=0$ for every $\alpha\in\hat R^+$, it follows that $w\lambda(c) = \lambda(c)$ for all $\lambda\in\hat P, w\in\widehat{\mathcal W}$. We also let
\begin{equation*}
c_{i,j}= \alpha_j(h_i), \quad i,j\in \hat I,
\end{equation*}
be the entries of the associated affine Cartan matrix. Recall also that, for $\lambda,\mu\in\hat P$, we have
\begin{equation*}
\mu\le\lambda \quad\text{if}\quad \lambda-\mu\in \hat Q^+ = \bigoplus_{i\in\hat I}\ \mathbb Z_{\ge 0}\ \alpha_i.
\end{equation*}

\subsection{Integrable Modules for Affine Algebras}\label{ss:intdem}
Given a $\hlie g$-module $V$ and $\mu\in\hlie h^*$, let $V_\mu$ denote the corresponding weight space of $V$:
\begin{equation*}
V_\mu = \{v\in V: hv=\mu(h)v \text{ for all } h\in\hlie h\}.
\end{equation*}
Given $\Lambda\in\hat P^+$, let $V(\Lambda)$ be an integrable irreducible $\hlie g$-module of highest weight $\Lambda$. Then
\begin{equation*}
\dim V(\Lambda)_{w\mu} = \dim V(\Lambda)_\mu \quad\text{for all}\quad w\in \widehat{\mathcal W}, \mu\in\hat P.
\end{equation*}
In particular, $\dim V(\Lambda)_{w\Lambda}=1$ for all $w\in \widehat{\mathcal W}$. The Demazure module $V_w(\Lambda)$ associated to $\Lambda\in\hat P^+$ and $w\in \widehat{\mathcal W}$ is the $\hlie b$-submodule of $V(\Lambda)$ generated by $V(\Lambda)_{w\Lambda}$. 

Let $\preccurlyeq$ denote the Bruhat order in $\widehat{\mathcal W}$. It is well-known (see \cite[Lemmas 1.3.20 and 8.3.3]{kumar}) that $(\widehat{\mathcal W},\preceq)$ is a directed poset,
\begin{gather}\notag 
w\preceq w' \quad\Rightarrow\quad V_w(\Lambda)\subseteq V_{w'}(\Lambda),\\ \label{e:BruhatDem}\text{and}\\ \notag
V(\Lambda) = \bigcup_{w\in \widehat{\mathcal W}} V_w(\Lambda).
\end{gather}
Since every element of $P$ is $\mathcal W$-conjugate to an element of $-P^+$, one easily sees that, for every $w\in\widehat{\mathcal W}$, there exists $w'\in{\mathcal W}$ such that
\begin{equation*}
\ell(w'w)=\ell(w)+\ell(w') \qquad\text{and}\qquad w'w\Lambda(h_i)\le 0 \quad\text{for all}\quad i\in I.
\end{equation*}
Thus, if we set
\begin{equation*}
{\widehat{\mathcal W}_\Lambda}^- = \{w\in \widehat{\mathcal W}: w\Lambda(h_i)\le 0 \text{ for all } i\in I\},
\end{equation*}
we get
\begin{equation}\label{e:demunion}
V(\Lambda) = \bigcup_{w\in {\widehat{\mathcal W}}_\Lambda^-} V_w(\Lambda).
\end{equation}
We refer to a Demazure module  of the form $V_w(\Lambda)$ for $w\in {\widehat{\mathcal W}}_\Lambda^-$ as a $\lie g$-stable Demazure module.
The reason for this terminology comes from the following well-known and easily established lemma.

\begin{lem}\label{l:gstdem}
	Let $\Lambda\in\hat P^+$ and $w\in \widehat{\mathcal W}$. The following are equivalent:
	\begin{enumerate}[(i)]
		\item $V_w(\Lambda)$ is a $\lie g[t]$-submodule of $V(\Lambda)$.
		\item $V_w(\Lambda)$ is a $\lie g$-submodule of $V(\Lambda)$.
		\item $\lie n^-V_w(\Lambda)_{w\Lambda} = 0$.
		\item $w\Lambda(h_i)\le 0$ for all $i\in I$.\endd
	\end{enumerate}
\end{lem}

Suppose $w\in {\widehat{\mathcal W}}_\Lambda^-$. Since $\hat P = \mathbb Z\Lambda_0\oplus P\oplus\mathbb C\delta$, there exist unique $\lambda\in P^+$ and $m\in\mathbb C$ such that
\begin{equation}\label{e:demnot}
w\Lambda = \ell\Lambda_0+w_0\lambda+m\delta,
\end{equation}
where $w_0$ is the longest element of $\mathcal W$ and $\ell$ is the level of $\Lambda$. Conversely, given $\ell\in\mathbb Z_{> 0},\lambda\in P^+$, and $m\in\mathbb C$, there exists $w\in \widehat{\mathcal W}$ and unique $\Lambda\in\hat P^+$ such that \eqref{e:demnot} holds. Since  $\sigma\delta=\delta$ for all $\sigma\in\widehat{\mathcal W}$, it follows that $V_w(\Lambda)$ and $V_w(\Lambda+z\delta)$ are isomorphic as (non-graded) $\lie g[t]$-modules for all $z\in\mathbb C$. Thus, for convenience, we assume that $\Lambda(d)\in\mathbb Z$ and, hence, $m\in\mathbb Z$. By \cite[Theorem 2]{CV}, we have an isomorphism of graded $\lie g[t]$-modules
\begin{equation}\label{e:cvdemrel}
V_w(\Lambda)\cong D(\ell,\lambda,m).
\end{equation}
It will be convenient to introduce the notation
\begin{equation}\label{e:Gamma}
\Gamma_\Lambda = \left\{(\lambda,m)\in P^+\times\mathbb Z:  \ell\Lambda_0+\lambda+m\delta \in\widehat{\mathcal W}\Lambda\right\}.
\end{equation}
Observe that
\begin{equation}\label{e:Gammatrans}
\Gamma_{\Lambda+s\delta} = \{(\lambda,r+s): (\lambda,r)\in\Gamma_\Lambda\}.
\end{equation}
Moreover, if $\Lambda'\in\hat P^+$ is a level-$\ell$ element distinct from $\Lambda$, then
\begin{equation}\label{e:capGamma}
  \Gamma_\Lambda\cap \Gamma_{\Lambda'} = \emptyset.
\end{equation}
Indeed, if $(\lambda,m)\in \Gamma_\Lambda\cap \Gamma_{\Lambda'}$, then $ \ell\Lambda_0+\lambda+m\delta$ would be in the ${\widehat{\mathcal W}}$-orbit of two distinct elements of $\hat P^+$, which is impossible.

\subsection{Demazure Flags}
We shall say that a $\lie g[t]$-module $V$ admits a Demazure flag if there exist $l>0, \lambda_j\in P^+, m_j\in\mathbb Z, j=1,\dots,l$, and a sequence of inclusions
\begin{equation}\label{e:demflag}
0 = V_0 \subset V_1 \subset \cdots \subset V_{l-1} \subset V_l = V \quad\text{with}\quad V_j/V_{j-1}\cong D(\ell_j, \lambda_j,m_j) \ \forall\ 1\le j\le l.
\end{equation}
If $\ell_j=\ell$ for some $\ell$ and all $j$, we say that such a sequence is a level-$\ell$ Demazure flag for $V$. Let $\mathbb V$ be a Demazure flag of $V$ as in \eqref{e:demflag} and, for a Demazure module $D$, define the multiplicity of $D$ in $\mathbb V$ by
\begin{equation*}
[\mathbb V:D] = \#\{1\le j\le l : V_j/V_{j-1}\cong D\}.
\end{equation*}
As observed in \cite[Lemma 2.1]{cssw:demflag}, the multiplicity does not depend on the choice of the flag and, hence, by abuse of language, we shift the notation from $[\mathbb V:D]$ to $[V:D]$. Also following \cite{cssw:demflag}, we consider the generating polynomial
\begin{equation*}
[V:D](q) = \sum_{m\in\mathbb Z}\ [V:\tau_mD]\ q^m \ \in\ \mathbb Z[q,q^{-1}].
\end{equation*}
In the category of non-graded $\lie g[t]$-modules we have $\tau_mD\cong D$ and, hence, one can also be interested in computing the ungraded multiplicity of $D$ in $V$ which is given by 
\begin{equation*}
[V:D](1) = \sum_{m\in\mathbb Z}\ [V:\tau_mD].
\end{equation*}

It  was proved in \cite[Theorem A]{naoi:weyldem} that $W(\lambda)$ has a level-$1$ Demazure flag for all $\lambda\in P^+$, thus generalizing \eqref{e:sliso}. The proof is based on results of \cite{joseph03,joseph06} which can be used to describe precisely the multiplicities of the Demazure modules in the flag. One of the main ingredients in the proof of \cite[Theorem A]{naoi:weyldem} is the following result (see \cite[Corollary 4.16]{naoi:weyldem}).

\begin{thm}\label{t:naoislf}
Suppose $\lie g$ is simply laced, let $\mu\in P^+$ and $\ell'>\ell\ge 0$. Then, $D(\ell,\mu)$ admits a level-$\ell'$ Demazure flag.\endd
\end{thm}

The proof of Theorem \ref{t:naoislf} in \cite{naoi:weyldem} uses a combinatorial description of certain Demazure crystals from \cite[Theorem 2.11]{joseph03} and 
a globalization result \cite[Corollary 5.10]{joseph06}. We recall the former in Section \ref{s:crystals}. As shown in \cite{naoi:weyldem}, once Theorem \ref{t:naoislf} is proved for $\lie g$ of type $A$, the existence of a level-$1$ Demazure flag for Weyl modules is then proved by applying the type $A$ case to the subalgebra generated by the simple short roots in an appropriate way. In particular, Theorem \ref{t:naoislf} is the only portion of the proof of \cite[Theorem A]{naoi:weyldem} depending on quantum groups arguments. 
A quantum-free proof of Theorem \ref{t:naoislf} for $\lie g=\lie{sl}_2$ was later given in \cite{cssw:demflag} which then leads to a quantum-free proof of the existence of level-$1$ Demazure flags for Weyl modules when $\lie g$ is of type $B$ or $G$. 

\begin{rem}
	As an application of the existence of a level-$1$ Demazure flag for $W(\lambda)$, it was proved in \cite{naoi:weyldem} that we have an isomorphism of $\lie g$-modules
	\begin{equation}\label{e:cpconj}
	W(\lambda)\cong \otm_{i\in I}^{} W(\omega_i)^{\otimes \lambda(h_i)},
	\end{equation}
	as conjectured in \cite{CPweyl} (for simply laced $\lie g$ this had been proved previously in \cite{chlo:wfd,foli:weyldem}). Both \cite[Theorem A]{naoi:weyldem} and \eqref{e:cpconj} were extended to the positive characteristic case in \cite{bmm,jm:hyper}. In particular, it follows from \cite[Theorem 1.5.2(b)]{bmm} that the multiplicities of Demazure modules in these flags are independent of the (algebraically closed) ground field.\hfill$\diamond$
\end{rem}

\subsection{Outer Multiplicities from Demazure Flags}
As mentioned in the introduction, if $\Upsilon,\Lambda\in\hat P^+$, then 
\begin{equation*}
V(\Upsilon)\otimes V(\Lambda) \cong \bigoplus_{\Phi\in\hat P^+} V(\Phi)^{\oplus m_\Phi} \quad\text{for some}\quad m_\Phi\in\mathbb Z_{\ge 0}.
\end{equation*}
Moreover, $m_\Phi\ne 0$ only if $\Phi\le \Upsilon+\Lambda$. Given a $\hlie g$-module $V$ such that
\begin{equation*}
V \cong \bigoplus_{\Phi\in\hat P^+} V(\Phi)^{\oplus m_\Phi} \quad\text{for some}\quad m_\Phi\in\mathbb Z_{\ge 0},
\end{equation*}
we set
\begin{equation*}
[V:V(\Phi)] = m_\Phi.
\end{equation*}
We are ready to state our main result.

\begin{thm}\label{t:multrel}
	If $\lie g$ is simply laced, $\Lambda,\Phi\in\hat P^+, V=V(\Lambda_0)\otimes V(\Lambda)$, and $\ell=\Lambda(c)$, then
	\begin{equation*}
	[V:V(\Phi)] =  \sum_{ (\lambda,r)\in\Gamma_\Phi}\  \max_{(\mu,s)\in\Gamma_\Lambda}\ [D(\ell,\mu,s): D(\ell+1,\lambda,r)].
	\end{equation*}
\end{thm}

Theorem \ref{t:multrel} will be proved in Section \ref{ss:multrel}. Given an affine Dynkin diagram automorphism $\sigma:\hat I\to\hat I$ and $j=\sigma(0)$, we can use Theorem \ref{t:multrel} to compute 
\begin{equation*}
[V(\Lambda_j)\otimes V(\Lambda):V(\Phi)] \qquad\text{for all}\qquad \Lambda,\Phi\in\hat P^+.
\end{equation*}
To do so, denote also by $\sigma$ the automorphism of $\hat P$ determined by
\begin{equation*}
\sigma(\Lambda_i)= \Lambda_{\sigma^{-1}(i)},\ i\in\hat I, \quad\text{and}\quad \sigma(\delta)=\delta.
\end{equation*}
Recall that the condition $\Phi\le\Lambda_j+\Lambda$ is equivalent to saying that
\begin{equation*}
\Lambda_j+\Lambda-\Phi =  \sum_{i\in\hat I} c_i\alpha_i \quad\text{for some}\quad c_i\in\mathbb Z_{\ge 0}
\end{equation*} 
and, in that case, the integers $c_i$ form the unique solution of the linear system
\begin{equation*}
c_0 = \Lambda(d)-\Phi(d), \qquad \sum_{i\in\hat I} c_{k,i} c_i = \Lambda(h_k)-\Phi(h_k)+\delta_{k,j} \qquad\text{for all}\qquad k\in\hat I.
\end{equation*}
Then,
\begin{equation}\label{e:multrelnot0}
[V(\Lambda_j)\otimes V(\Lambda):V(\Phi)] = [V(\Lambda_0)\otimes V(\sigma(\Lambda)):V(\sigma(\Phi) +(c_0-c_j)\delta)]. 
\end{equation}
In particular, if $\lie g$ is of type $A$, this gives a way of computing 
$[V(\Lambda_i)\otimes V(\Lambda):V(\Phi)]$ for all $i\in\hat I,\Lambda,\Phi\in\hat P^+$.  Although the proof of \eqref{e:multrelnot0} utilizes standard arguments with Dynkin diagram automorphisms, we include it in Section \ref{ss:multrelnot0} for the reader's convenience. 

\subsection{The $\lie{sl}_2$-Case}

As of this moment, the only case for which the multiplicities in $\lie g$-stable Demazure flags are described in a more explicit manner is for $\lie g=\lie{sl}_2$ (see Remark \ref{r:main} below for more comments). In that case, we are able to obtain a more explicit expression for the right-hand side of the formula given in Theorem \ref{t:multrel}. Thus, we assume $\lie g=\lie{sl}_2$, $\hat I=\{0,1\}$, and, when convenient, we identify $P$ with $\mathbb Z$ by sending $\omega_1$ to $1$. It follows from the results of \cite{cssw:demflag} that, for $\ell'\ge \ell$, we have
\begin{equation}\label{e:parity}
[D(\ell,\lambda) : D(\ell',\lambda)](q) = 1 \quad\text{and}\quad [D(\ell,\lambda) : D(\ell', \mu)](q) = 0 \quad\text{if}\quad \lambda-\mu \notin 2\mathbb Z_{\ge 0}.
\end{equation}
Consider the generating series
\begin{equation*}
A_\lambda^{\ell,\ell'}(x,q) = \sum_{m\ge 0} [D(\ell,\lambda+2m) : D(\ell', \lambda)](q)\ x^m =  
\sum_{m\ge 0}\sum_{r\in\mathbb Z} [D(\ell,\lambda+2m) : D(\ell', \lambda,r)]\ q^r x^m.
\end{equation*}
For $1\le\ell\le\ell'\le 3$, it was shown in \cite{bcsv:dfcpmt} (see Theorems 1.5 and 1.6) that $A_\lambda^{\ell,\ell'}(x,q)$ are related to partial and mock theta functions. In particular, $A_\lambda^{\ell,\ell'}(x,1)$ is a generating series for ungraded multiplicities and it was shown in \cite[Section 1.3]{bcsv:dfcpmt} that it can be expressed in terms of Chebyshev polynomials. To make the connection with Theorem \ref{t:multrel}, let $\ell'=\ell+1$ and write $A_\lambda^\ell(x,q)$ instead of $A_\lambda^{\ell,\ell+1}(x,q)$. Set
\begin{equation}\label{e:defalfa}
\alpha_\lambda^\ell(m,r) = [D(\ell,\lambda+2m) : D(\ell+1, \lambda,r)]
\end{equation}
for all $\lambda,\ell,r\in\mathbb Z, \lambda\ge 0,\ell> 0, m\in\mathbb Q$. In particular,
\begin{equation*}
\alpha_\lambda^\ell(m,r) = 0 \quad\text{if}\quad m\notin\mathbb Z_{\ge 0}
\end{equation*}
and
\begin{equation*}
A_\lambda^\ell(x,q) = \sum_{m\ge 0}\sum_{r\in\mathbb Z}\ \alpha_\lambda^\ell(m,r)\ x^m\ q^r.
\end{equation*}

\begin{prop}\label{c:sl2l}
	Let $\Lambda\in\hat P^+, i\in\hat I$, and $\sigma:\hat I\to\hat I$ be the nontrivial diagram automorphism. Then, for all $\Phi\in\hat P^+$, we have
	\begin{equation*}
	[V(\Lambda_i)\otimes V(\Lambda):V(\Phi)] =  \sum_{ (\lambda,r)\in\Gamma_{\sigma^i(\Phi)}}\ \lim_{k\to\infty}\  \alpha_\lambda^\ell \left(k\ell+\frac{m_{i}-\lambda}{2}\ ,\ r_i+ r+ k(k\ell+m_i)-s\right).
	\end{equation*}	
	where $\ell=\Lambda(c), m_i=\Lambda(h_{1-i}), r_i =i\ \frac{\Lambda(h_0)-\Phi(h_0)}{2}$, and $s=\Lambda(d)$.
\end{prop}

The proof of Proposition \ref{c:sl2l} will be given in Section \ref{ss:sl2l}. It is an almost immediate corollary of Theorem \ref{t:multrel} and Weyl group orbits computations performed in Section \ref{ss:sl2orbits}.

In the case that $\ell=1$, we can explicitly compute the limits appearing in Proposition \ref{c:sl2l} in terms of partitions with bounded parts. Given $m>0,k\ge 0$, denote by $\bro_k(m)$ the number of partitions of $m$ with all parts bounded by $k$. Set also $\bro_k(0)=1$ and $\bro_k(m)=0$ if $m<0$.
For a nonnegative real number $a$ we let $\lfloor a\rfloor$ be the integer part of $a$ and $\lceil a\rceil=\lfloor a\rfloor+1$.

\begin{prop}\label{t:sl21}
	Let $V=V(\Lambda_0)\otimes V(\Lambda_i), i\in\hat I$, and $\Phi\in\hat P^+$. Then, $[V:V(\Phi)]\ne 0$ only if $\Phi=2\Lambda_0+(2j+i)\omega_1-s\delta$ for some $0\le j\le \delta_{i,0},s\ge j$. In that case, for $i=0$,
	\begin{equation*}
	[V:V(\Phi)] = \sum_{l=0}^{\lfloor L\rfloor}\ \bro_{2l+j}(s-2l^2-2jl-j) \qquad\text{where}\qquad L=\frac{-j+\sqrt{2s-j}}{2},
	\end{equation*}	
	and, for $i=1$, 
	\begin{align*}
	[V:V(\Phi)]  =  \sum_{l=0}^{\lfloor L\rfloor}  \bro_{2l}(s-2l^2+l) \qquad\text{where}\qquad L=\frac{1+\sqrt{8s+1}}{4}.
	\end{align*}
\end{prop}

Proposition \ref{t:sl21} will be proved in Section \ref{ss:sl21}. Combining \eqref{e:multrelnot0} with Proposition \ref{t:sl21} for the case $i=0$, one easily obtains the following corollary.

\begin{cor}
	Let $V=V(\Lambda_1)\otimes V(\Lambda_1)$ and $\Phi\in\hat P^+$. Then, $[V:V(\Phi)]\ne 0$ only if $\Phi=2\Lambda_{1-j}-s\delta$ for some $j\in\hat I,s\ge 0$, and, in that case,
	\begin{equation*}
	[V:V(\Phi)] = \sum_{l=0}^{\lfloor L\rfloor}\ \bro_{2l+j}(s-2l^2-2jl) \qquad\text{where}\qquad L=\frac{-j+\sqrt{2s+j}}{2}.
	\end{equation*}	\endd
\end{cor}

In \cite[Theorem 2.1]{mw:tp2}, the numbers $[V:V(\Phi)]$ with $\lie g$ of type $A$ and $V=V(\Lambda_0)\otimes V(\Lambda_i), i\in\hat I$, were also expressed as the cardinality of a set of partitions satisfying certain conditions which do not coincide with the ones in Proposition \ref{t:sl21} (see \cite[Lemma 2.3]{mw:tp2}). For instance, setting 
\begin{equation*}
\Phi^{i}_{j,s} = 2\Lambda_0+(2j+i)\omega_1-s\delta \quad\text{for}\quad i\in\hat I,\ 0\le j\le \delta_{i,0},\ s\ge j,
\end{equation*}
it follows from \cite{mw:tp,mw:tp2} that
\begin{equation}\label{e:mw00}
[V(\Lambda_0)\otimes V(\Lambda_i):V(\Phi^i_{j,s})] = \bro_{\ne}^i(2s-j),
\end{equation}
where $\bro_{\ne}^i(k)$ is the cardinality of the set of partitions of $k$ with all parts distinct and having parity different from that of $i$.
Combining this with Proposition \ref{t:sl21}, we get the following partition identity:
\begin{equation}\label{e:partrel}
\bro_{\ne}^i(2s-j) = \sum_{l\ge 0} \bro_{2l+j}(s-2l^2-(2j-i)l-j) \quad\text{for}\quad i\in\hat I,\ 0\le j\le \delta_{i,0},\ s\ge j.
\end{equation}
Such results can also be used to obtain $q$-identities. To do that, following \cite{mw:tp,mw:tp2}, set
\begin{equation}
b^i_{j,s} = [V(\Lambda_0)\otimes V(\Lambda_i): V(\Phi^i_{j,s})] \qquad\text{and}\qquad B^i_j(q) = \sum_{s=|i-j|}^\infty b^i_{j,s}\ q^{s-|i-j|}.
\end{equation}
Using \eqref{e:mw00} and \cite[Theorem 1.1]{andrews}, it was shown in \cite{mw:tp} that
\begin{equation}\label{e:Bformula}
B_j^0(q) = \sum_{l\ge 0} \frac{q^{2l(l+1-j)}}{(q)_{2l+1-j}} 
\qquad\text{where}\qquad (a)_k = \prod_{r=0}^{k-1}(1-aq^r).
\end{equation}
If one starts from Proposition \ref{t:sl21} instead of \eqref{e:mw00} and applies \cite[Theorem 1.1]{andrews} in the same manner, one gets a new proof of \eqref{e:Bformula}.  As pointed out in \cite{mw:tp}, combining \eqref{e:Bformula} with the results from \cite{fein} leads to new proofs of certain identities listed in \cite{slater}.

\begin{rem}\label{r:main}
	The two main tools, besides Theorem \ref{t:multrel}, used in the proof of Proposition \ref{t:sl21} are an explicit description of the set $\Gamma_\Phi$ (given in Section \ref{ss:sl2orbits} for $\Phi$ of any level) and the fact, proved in \cite{cssw:demflag}, that the generating function $[D(1,\mu),D(2,\lambda)](q)$ is essentially a Gaussian binomial (see \eqref{e:csswmf2} below). The generating function $[D(2,\mu),D(3,\lambda)](q)$ has been expressed as a finite sum of products of two Gaussian binomials in \cite[Proposition 1.4]{bcsv:dfcpmt}. Therefore, following the steps of the proof of Proposition \ref{t:sl21}, one can obtain a formula for tensor products of the form $V(\Lambda_i)\otimes V(\Lambda)$ with $\Lambda$ of level $2$ in terms of partitions with bounded parts. For higher-level $\Lambda$, it is unclear to us at this moment if the information on $[D(\ell,\mu),D(\ell+1,\lambda)](q)$ given in \cite{bcsv:dfcpmt,cssw:demflag} is sufficient to actually compute the limit in Proposition \ref{c:sl2l}.
	
	For $\lie g=\lie{sl}_3$, $[D(1,\mu),D(2,\lambda)](q)$ has been expressed as a finite sum of products of two Gaussian binomials in \cite[Theorem 16]{wand}. Thus, once a description of $\Gamma_\Phi$ for $\Phi$ of level $2$ is given, the steps of the proof of Proposition \ref{t:sl21} also lead to 
	a formula for outer multiplicities of tensor products of the form $V(\Lambda_i)\otimes V(\Lambda_j)$ in terms of partitions with bounded parts.  As mentioned above, such outer multiplicities were computed in \cite{mw:tp2} using very different methods. Moreover, we have seen that, even in the $\lie{sl}_2$ case, the two methods lead to different expressions in terms of partitions. Thus, more intricate partition identities than \eqref{e:partrel} are expected. We shall perform such computations and comparisons with the results of \cite{mw:tp2} in a future publication.  \endd
\end{rem}

\section{Crystal Approach to Demazure Flags}\label{s:crystals}

\subsection{The PMK Theorem} We now review a crucial piece of the background to be used in the proof of Theorem \ref{t:multrel}.

The notion of Demazure flags  can be extended from $\lie g$-stable Demazure modules to general ones in the obvious way. Note that the highest-weight space of $V(\Lambda)$ is a Demazure module: $V(\Lambda)_\Lambda = V_{id}(\Lambda)$. For historical comments and precise references for the proof of the next theorem, see \cite{joseph03,joseph06,mat:pos}. The theorem is referred to as PMK Theorem in \cite{joseph03} after independent works of P. Polo, O. Mathieu, and W. van der Kallen.

\begin{thm}\label{t:PMK}
For any $\Upsilon,\Lambda\in\hat P^+$ and $w\in \widehat{\mathcal W}$, $V_{id}(\Upsilon)\otimes V_w(\Lambda)$ admits a Demazure flag.\endd
\end{thm}

Let $\mu,\ell$ be as in Theorem \ref{t:naoislf} and let $\Lambda\in\hat P^+, w\in \widehat{\mathcal W}$ be such that
\begin{equation*}
w\Lambda = w_0\mu + \ell\Lambda_0.
\end{equation*}
The proof of Theorem \ref{t:naoislf} with $\ell'=\ell+1$ uses a global basis version of  Theorem \ref{t:PMK} (\cite[Corollary 5.10]{joseph06}) with $\Upsilon=\Lambda_0$ (for general $\ell'$ one then iterates the procedure). More precisely, let
\begin{equation*}
0= Y_0\subset Y_1 \subset\cdots\subset Y_l =  V_{id}(\Lambda_0)\otimes V_w(\Lambda)
\end{equation*}
be a Demazure flag, say,
\begin{equation}\label{e:mainsl}
Y_j/Y_{j-1} \cong V_{w_j}(\Phi_j) \quad\text{for}\quad j=1,\dots,l.
\end{equation}
Then, the proof of Theorem \ref{t:naoislf} implies that there exists a Demazure flag
\begin{equation*}
0= W_0\subset W_1 \subset\cdots\subset W_l =  V_w(\Lambda) \cong D(\ell,\mu)
\end{equation*}
such that
\begin{equation}
W_j/W_{j-1} \cong Y_j/Y_{j-1} \quad\text{and, moreover,}\quad \Phi_j(c)=\ell+1 \quad\text{for all}\quad j=1,\dots,l.
\end{equation}
Combinatorial versions of Theorem \ref{t:PMK} were proved in \cite[Theorem 2.11]{joseph03}, \cite[Section 2.4]{llm}, \cite[Theorem 9(b)]{lit:qFrsm} from which it is, in principle, possible to compute the elements $w_j$ and $\Phi_j$ of \eqref{e:mainsl}. We will use the version from \cite{joseph03} (see  Theorem \ref{t:cPMK}  and \eqref{e:cmainsl} below) in the proof of Theorem \ref{t:multrel}. The underlying combinatorial context is that of crystals which we briefly review next.

\subsection{Crystals} We begin by recalling some basics on crystal theory and refer to \cite{jos:book,kas:CB,kas:dem} and references therein for further details. A crystal (associated to the Cartan data of $\hlie g$) is a nonempty set $B$ with maps $\wt:B\to\hat P, \varepsilon_i,\varphi_i:B\to\mathbb Z\cup\{-\infty\}, e_i,f_i:B\cup\{0\}\to B\cup\{0\}, i\in\hat I$, satisfying, in particular,  the following properties:
\begin{enumerate}[(1)]
\item $e_i0=f_i0=0$;
\item given $i\in\hat I, b,b'\in B$, one has $e_ib = b' \Leftrightarrow f_ib'=b$;
\item if $e_ib\ne 0$ for some $i\in\hat I, b\in B$, then $\wt(e_ib) = \wt(b)+\alpha_i$;
\item $\varphi_i(b)=\varepsilon_i(b)+\wt(b)(h_i)$ for all $b\in B$.
\end{enumerate}
A crystal $B$ is said to be upper normal if
\begin{enumerate}[(5)]
\item $\varepsilon_i(b) = \max\{k: e_i^kb\ne 0\}$ for all $i\in\hat I, b\in B$.
\end{enumerate}
Similarly one defines lower normal crystal by replacing $\varepsilon_i$ by $\varphi_i$ and $e_i$ by $f_i$ above. If $B$ is both upper and lower normal, it is simply said that $B$ is a normal crystal. Given $\mu\in\hat P$, set $B_\mu=\{b\in B:\wt(b)=\mu\}$. We shall always assume that $\# B_\mu<\infty$ for all $\mu\in\hat P$. If $B$ is normal,
\begin{equation*}
\# B_\mu = \# B_{w\mu} \quad\text{for all}\quad \mu\in\hat P.
\end{equation*}
Let $\mathcal E\subseteq {\rm End_{Set}}(B\cup\{0\})$ be the submnoid generated by $e_i,i\in\hat I$, and similarly define $\mathcal F$. An element $b\in B$ is said to be primitive if $e_ib=0$ for all $i\in\hat I$ or, equivalently, $\mathcal Eb=\{0\}$. $B$ is said to be a highest-weight crystal of highest weight $\lambda\in\hat P$ if there exists a primitive element $b\in B_\lambda$ such that $B=\mathcal Fb\setminus\{0\}$. In that case, for simplicity, we shall write $B=\mathcal Fb$.

Given a family of crystals $B_j, j\in J,$ we equally use the notation $\wt,e_i,f_i,\varepsilon_i,\varphi_i,\mathcal E,\mathcal F$ for each of the $B_j$. This will not create confusion. A morphism from the crystal $B_1$ to the crystal $B_2$ is a map $\psi:B_1\cup\{0\}\to B_2\cup\{0\}$ such that
\begin{enumerate}[(1)]
\item $\psi(0)=0$;
\item $\wt(\psi(b))=\wt(b)$ for all $b\in B_1$ and similarly for $\varepsilon_i,\varphi_i$ in place of $\wt$ for all $i\in\hat I$;
\item $\psi(e_ib) = e_i\psi(b)$ if $e_ib\ne 0$ and similarly for $f_i$ in place of $e_i$.
\end{enumerate}
If $\psi$ is injective, $B_1$ is said to be a subcrystal of $B_2$. In the case that $\psi$ commutes with $e_i$ and $f_i$ for all $i\in\hat I$, $\psi$ is said to be a strict morphism and $B_1$ is said to be a strict subcrystal of $B_2$. An isomorphism of crystals is a bijective morphism.

The character of a crystal is defined in the obvious way. For each $\Lambda\in\hat P^+$, there exists a unique, up to isomorphism, normal crystal $B(\Lambda)$ of highest weight $\Lambda$. Moreover,
\begin{equation*}
\# B(\Lambda)_\mu = \dim V(\Lambda)_\mu \quad\text{for all}\quad \mu\in\hat P.
\end{equation*}

Given two crystals $B_1$ and $B_2$, the set $B_1\times B_2$ can be equipped with the structure of a crystal, denoted $B_1\otimes B_2$. However, we will not need the precise definition of the structure here. 
For $\lambda\in\hat P^+$, let $b_\lambda$ be the highest-weight generator of $B(\lambda)$ and, for $\mu\in\hat P^+$, set 
$$B(\mu)^\lambda = \{b\in B(\mu): b_\lambda\otimes b \text{ is primitive}\}.$$ 

\begin{thm}\label{t:tpcrystal}
Let $\lambda,\mu\in\hat P^+$. Then, every primitive element of $B(\lambda)\otimes B(\mu)$ is of the form $b_\lambda\otimes b$ with $b\in B(\mu)$. Moreover, if  $b\in B(\mu)^\lambda$, then $\wt(b_\lambda\otimes b)\in\hat P^+$, $\mathcal F(b_\lambda\otimes b)$ is a strict subcrystal isomorphic to $B(\lambda+\wt(b))$, and
\begin{equation*}
B(\lambda)\otimes B(\mu) = \bigsqcup_{b\in B(\mu)^\lambda} \mathcal F(b_\lambda\otimes b).
\end{equation*}
\end{thm}

\subsection{Demazure Crystals}\label{ss:demc} Given a crystal $B$ and a reduced expression $w=s_{i_1}\cdots s_{i_l}$ of $w\in\widehat{\mathcal W}$, consider
\begin{equation*}
\mathcal F_w =  \{f_{i_1}^{k_1}\cdots f_{i_l}^{k_l}: k_j\in\mathbb Z_{\ge 0}, j=1,\dots,l\}\subseteq\mathcal F.
\end{equation*}
In principle, $\mathcal F_w$ depends on the choice of reduced expression for $w$. However, if $B$ is the union of normal highest-weight subcrystals, then it is known that $\mathcal F_w$ is well-defined (see \cite{jos:book,joseph03,kas:dem} and references therein). The Demazure crystal associated to $w\in\widehat{\mathcal W}$ and $\Lambda\in P^+$ is the subset of $B(\Lambda)$ given by
\begin{equation*}
B_w(\Lambda) = \mathcal F_w b_\Lambda.
\end{equation*}
It is not a subcrystal of $B(\Lambda)$ in general, but it is $\mathcal E$-stable, i.e., $$\mathcal E B_w(\Lambda)\subseteq B_w(\Lambda)\cup\{0\}.$$
It was proved in \cite[Section 4.6]{joseph03} that 
\begin{equation*}
\ch(V_w(\Lambda)) = \ch(B_w(\Lambda)).
\end{equation*}
Moreover, we also have
\begin{equation}\label{e:demunionc}
w\preceq w' \quad\Rightarrow\quad B_w(\Lambda)\subseteq B_{w'}(\Lambda) \qquad\text{and}\qquad B(\Lambda) = \bigcup_{w\in {\widehat{\mathcal W}}_\Lambda^-} B_w(\Lambda).
\end{equation}

The notion of morphisms of Demazure crystals is defined in the obvious way.
Set
\begin{equation*}
B_w(\mu)^\lambda = B(\mu)^\lambda\cap B_w(\mu), \quad \lambda,\mu\in\hat P^+,\ w\in\widehat{\mathcal W}.
\end{equation*}
The following is the combinatorial version of Theorem \ref{t:PMK} given in \cite[Theorem 2.11]{joseph03}.

\begin{thm}\label{t:cPMK}
For every $\Upsilon,\Lambda\in\hat P^+$, and $w\in\widehat{\mathcal W}$, the subset $b_\Upsilon\otimes B_w(\Lambda)$ of $B(\Upsilon)\otimes B(\Lambda)$ is a disjoint union of Demazure crystals. More precisely, for each $b\in B_w(\Lambda)^\Upsilon$, there exists $w_b\in\widehat{\mathcal W}$ such that
\begin{equation*}
b_\Upsilon\otimes B_w(\Lambda) = \bigsqcup_{b\in B_w(\Lambda)^\Upsilon} \mathcal F_{w_b}(b_\Upsilon\otimes b) \qquad\text{and}\qquad
\mathcal F_{w_b}(b_\Upsilon\otimes b) \cong B_{w_b}(\Upsilon+\wt(b)).
\end{equation*}\endd
\end{thm}

In the language of \eqref{e:mainsl}, we have $l=\#B_w(\Lambda)^\Upsilon$ and, for each $j=1,\dots,l$, there exists unique $b\in B_w(\Lambda)^\Upsilon$ such that \begin{equation}\label{e:cmainsl}
w_j = w_b \qquad\text{and}\qquad \Phi_j = \Upsilon+\wt(b).
\end{equation}

\subsection{The Proof of Theorem \ref{t:multrel}}\label{ss:multrel}
	It follows from Theorem \ref{t:tpcrystal} that 
	\begin{equation*}
	[V:V(\Phi)] = \# B(\Lambda)^{\Lambda_0}_\Phi \quad\text{where}\quad B(\Lambda)^{\Lambda_0}_\Phi = \{b\in B(\Lambda)^{\Lambda_0}: \Phi = \wt(b)+{\Lambda_0}\}.
	\end{equation*}
		
	On the other hand, since $(\widehat{\mathcal W},\preceq)$ is a directed poset, for every $w,w'\in\hat{\mathcal W}$, there exists $w''\in\widehat{\mathcal W}$ such that $w\preceq w''$ and $w'\preceq w''$. Moreover, we can assume $w''\in\widehat{\mathcal W}_\Lambda^-$. It then easily follows from \eqref{e:demunionc} that there exists $w\in{\widehat{\mathcal W}}_\Lambda^-$ such that 
	\begin{equation}\label{e:gotomax}
	B(\Lambda)^{\Lambda_0}_\Phi\subseteq B_w(\Lambda)^{\Lambda_0}.
	\end{equation}
	For each $b\in B(\Lambda)^{\Lambda_0}_\Phi$, let $w_b$ be the element given by Theorem \ref{t:cPMK}. Thus,
	\begin{equation*}
	\mathcal F_{w_b}(b_{\Lambda_0}\otimes b) \cong B_{w_b}(\Phi).
	\end{equation*} 
	It follows from the proof of Theorem \ref{t:naoislf} (see comments around \eqref{e:mainsl}) that, for all $b\in B(\Lambda)^{\Lambda_0}_\Phi$,
	\begin{equation*}
	V_{w_b}(\Phi)\cong D(\ell+1,\lambda_b,r_b) \quad\text{for some}\quad \lambda_b\in P^+, r_b\in\mathbb Z.
	\end{equation*}
	Moreover, 
	\begin{equation}
	[V_w(\Lambda):D(\ell+1,\lambda,r)] = \#B(\Lambda)^{\Lambda_0}_\Phi(\lambda,r)
	\end{equation}
	where, for $(\lambda,r)\in\Gamma_\Phi$, 
	\begin{equation*}
	B(\Lambda)^{\Lambda_0}_\Phi(\lambda,r) = \{b\in B(\Lambda)^{\Lambda_0}_\Phi: w_b\Phi = (\ell+1)\Lambda_0+w_0\lambda+r\delta\}.
	\end{equation*}
	Therefore,
	\begin{equation*}
	[V:V(\Phi)] = \#B(\Lambda)^{\Lambda_0}_\Phi = \sum_{ (\lambda,r)\in\Gamma_\Phi}\#B(\Lambda)^{\Lambda_0}_\Phi(\lambda,r) =  \sum_{ (\lambda,r)\in\Gamma_\Phi} [V_w(\Lambda): D(\ell+1,\lambda,r)].
	\end{equation*}
	 We claim that, given $(\lambda,r)\in\Gamma_\Phi$, 
	\begin{equation}\label{e:wmax}
	[V_w(\Lambda): D(\ell+1,\lambda,r)] \ge [V_{w'}(\Lambda): D(\ell+1,\lambda,r)] \quad\text{for all}\quad w'\in\widehat{\mathcal W}^-_\Lambda.
	\end{equation}
	This implies
	\begin{equation*}
	[V:V(\Phi)] = \sum_{ (\lambda,r)\in\Gamma_\Phi}\  \max_{w\in{\widehat{\mathcal W}}_\Lambda^-}\ [V_w(\Lambda): D(\ell+1,\lambda,r)],
	\end{equation*}	
	which is clearly equivalent to the expression given in the statement of the theorem. Thus, it remains to prove \eqref{e:wmax}.

	Fix $(\lambda,r)\in\Gamma_\Phi$. If $w'\prec w$, \eqref{e:wmax} is clear from the first part of \eqref{e:BruhatDem}. If $w\prec w'$, recall that $B_w(\Lambda)\subseteq B_{w'}(\Lambda)$ by \eqref{e:demunionc}. It then follows from \eqref{e:gotomax} that
	\begin{equation*}
	 \wt(b)+\Lambda_0\ne \Phi \qquad\text{for all}\qquad b\in B_{w'}(\Lambda)^{\Lambda_0}\setminus B_{w}(\Lambda)^{\Lambda_0}. 
	\end{equation*}
	Setting $\Phi_b=\wt(b)+\Lambda_0$, it follows from \eqref{e:capGamma} that
	\begin{equation*}
	 \Gamma_\Phi\cap  \Gamma_{\Phi_b} = \emptyset \qquad\text{for all}\qquad b\in B_{w'}(\Lambda)^{\Lambda_0}\setminus B_{w}(\Lambda)^{\Lambda_0}.
	\end{equation*}
	This, together with \eqref{e:isodeml}, implies 
	\begin{equation*}
	  D(\ell+1,\mu,s) \ncong D(\ell+1,\lambda,r) \qquad\text{for all}\qquad  (\mu,s)\in\Gamma_{\Phi_b},\ b\in B_{w'}(\Lambda)^{\Lambda_0}\setminus B_{w}(\Lambda)^{\Lambda_0}.
	\end{equation*}
	Hence, equality holds in \eqref{e:wmax} in this case. Finally, if $w$ and $w'$ are not comparable, let $w''\in\widehat{\mathcal W}_\Lambda^-$ be such that $w\preceq w''$ and $w'\preceq w''$. The cases of \eqref{e:wmax}  which were already proved imply
	\begin{gather*}
	  [V_w(\Lambda): D(\ell+1,\lambda,r)] = [V_{w''}(\Lambda): D(\ell+1,\lambda,r)]  \\  \qquad\text{and}\qquad \\
	  [V_{w'}(\Lambda): D(\ell+1,\lambda,r)] \le [V_{w''}(\Lambda): D(\ell+1,\lambda,r)],
	\end{gather*}
	completing the proof of \eqref{e:wmax}.

\subsection{Proof of \eqref{e:multrelnot0}}\label{ss:multrelnot0} Let $x_0^\pm = x_\theta^\mp\otimes t^{\pm 1}$ and recall that the derived algebra $\hlie g'$ of $\hlie g$ is generated by $x_i^\pm, i\in\hat I$, and that $\hlie g=\hlie g'\oplus\mathbb Cd$.
For a Dynkin diagram automorphism $\varsigma:\hat I\to\hat I$, consider the corresponding automorphism of $\hlie g'$, also denoted by $\varsigma$, which is determined by the assignments
\begin{equation*}
x^\pm_i\mapsto x^\pm_{\varsigma(i)}, \quad i\in\hat I.
\end{equation*}
In particular, 
\begin{equation*}
h=\sum_{i\in\hat I} a_ih_i \quad\Rightarrow\quad \varsigma(h) = \sum_{i\in\hat I} a_ih_{\varsigma(i)}.
\end{equation*}
Given an integrable $\hlie g$-module $V$, let $V^\varsigma$ be the $\hlie g'$-module obtained from $V$ by pulling-back the action by $\varsigma$. We shall denote by $v^\varsigma$ the element $v\in V$ when regarded as an element of $V^\varsigma$. Thus, 
\begin{equation*}
xv^\varsigma = \varsigma(x)v \quad\text{for all}\quad x\in\hlie g',\ v\in V.
\end{equation*}
One easily checks that, for all $\mu\in\hat P$, we have
\begin{equation}\label{e:h'action}
v\in V_\mu \quad\Rightarrow\quad hv^\varsigma = (\varsigma(\mu)(h))v^\varsigma \quad\text{for all}\quad h\in\hlie h',
\end{equation}
where $\hlie h'=\hlie h\cap\hlie g'$. It follows that we have an isomorphism of $\hlie g'$-modules
\begin{equation*}
V(\Lambda)^\varsigma \cong V(\varsigma(\Lambda)) \quad\text{for all}\quad \Lambda\in\hat P^+.
\end{equation*}
Moreover, it is also clear that
\begin{equation*}
(V\otimes W)^\varsigma \cong V^\varsigma\otimes W^\varsigma.
\end{equation*}

Let $\sigma,j$, and $\Lambda$ be as in \eqref{e:multrelnot0} and consider the above discussion with $\varsigma=\sigma^{-1}$ which implies that we have an isomorphism of $\hlie g'$-modules
\begin{equation*}
\psi: V(\Lambda_j)\otimes V(\Lambda)\to (V(\Lambda_0)\otimes V(\sigma(\Lambda)))^{\sigma^{-1}}.
\end{equation*}
To shorten notation, set 
\begin{equation*}
V= V(\Lambda_0)\otimes V(\sigma(\Lambda)) \qquad\text{and}\qquad W = V(\Lambda_j)\otimes V(\Lambda).
\end{equation*}
We can use $\psi$ to turn $V^{\sigma^{-1}}$ into a $\hlie g$-module by defining its graded components to be the image of those of $W$, thus making $\psi$ be an isomorphism of $\hlie g$-modules. For each weight $\Phi$ of $W$, let $\Phi^\sigma$ be the unique weight of $V$ such that
\begin{equation*}
\psi(W_\Phi) = V_{\Phi^\sigma}.
\end{equation*}
Note that the definition of $\sigma$ on $\delta$ has been chosen so that 
\begin{equation*}
\Phi = \Lambda_j+\Lambda \quad\Rightarrow\quad \Phi^\sigma = \Lambda_0+\sigma(\Lambda).
\end{equation*}

It is now clear that
\begin{equation*}
[W:V(\Phi)] = [V:V(\Phi^\sigma)]
\end{equation*}
and we are left to show 
\begin{equation*}
\Phi^\sigma = \sigma(\Phi) + (c_0-c_j)\delta
\end{equation*}
where the numbers $c_j$ are defined just before \eqref{e:multrelnot0}. Since \eqref{e:h'action} implies 
\begin{equation*}
\Phi^\sigma(h) = \sigma(\Phi)(h) \quad\text{for all}\quad h\in\hlie h',
\end{equation*}
it remains to show 
\begin{equation*}
\Phi^\sigma(d) = \sigma(\Phi)(d) +(c_0-c_j).
\end{equation*}
Let $v$ be a nonzero vector in the top weight space of $V$. Then,
\begin{equation*}
x_{i_1}^-\cdots x_{i_l}^- v^{\sigma^{-1}} = x_{\sigma^{-1}(i_1)}^-\cdots x_{\sigma^{-1}(i_l)}^- v
\end{equation*}
which proves
\begin{equation*}
\Phi^\sigma = \Lambda_0+\sigma(\Lambda) - \sum_{i\in\hat I} c_i\alpha_{\sigma^{-1}(i)}.
\end{equation*}
Since 
\begin{equation*}
\sigma(\Phi) = \Lambda_0+\sigma(\Lambda) - \sum_{i\in\hat I}c_i\sigma(\alpha_i)
\end{equation*}
and a straightforward computation shows that
\begin{equation*}
\sigma(\alpha_i) = \alpha_{\sigma^{-1}(i)} +(\delta_{i,0}- \delta_{i,j}) \delta,
\end{equation*}
the proof is complete.

\section{The $\lie{sl}_2$-Case}\label{s:sl2}

\subsection{Weyl Group Orbits}\label{ss:sl2orbits}
Given $m\in\mathbb Z_{\ge 0}$, we simplify notation and write $W(m)$ instead of $W(m\omega_1)$ and similarly for $D(\ell,m,r)$. 

Note that
\begin{equation*}
\omega_1 = \Lambda_1 - \Lambda_0, \quad \theta=\alpha_1 = 2\omega_1 = 2\Lambda_1 - 2\Lambda_0, \quad \alpha_0 = \delta-2\omega_1 = \delta+2\Lambda_0-2\Lambda_1,
\end{equation*}
\begin{equation*}
s_0\Lambda_0 = \Lambda_0 +2\omega_1-\delta = -\Lambda_0 + 2\Lambda_1 -\delta, \qquad s_1\Lambda_1 = \Lambda_0-\omega_1 = 2\Lambda_0-\Lambda_1,
\end{equation*}
and $s_i\Lambda_j = \Lambda_j$ for $i\ne j$. Given $k\in\mathbb Z$, set
\begin{equation*}
\sigma_{k} = (s_1s_0)^ks_1
\end{equation*}
and recall that
\begin{equation*}
\widehat{\mathcal W} = \{\sigma_k, \sigma_ks_1 :k\in\mathbb Z\}.
\end{equation*}
The $\widehat{\mathcal W}$-orbit of $\Lambda_i,i\in\hat I$, is then explicitly described:
\begin{equation}\label{e:orbitLambdai2}
\sigma_k\Lambda_i=\sigma_{k+i}s_1\Lambda_i  = \Lambda_0-(2k+i)\omega_1-k(k+i)\delta \quad\text{for all}\quad k\in\mathbb Z.
\end{equation}
In particular, the $\lie g$-stable Demazure modules inside $V(\Lambda_i)$ are
\begin{equation}\label{e:weylsasdemsl2}
V_{\sigma_k}(\Lambda_i)\cong D(1,2k+i,-k(k+i)) \cong \tau_{-k(k+i)}W(2k+i) \quad\text{for}\quad k\ge 0.
\end{equation}

More generally, recall that every $\Lambda\in\hat P^+$ can be written as  $\Lambda = \ell\Lambda_0+m\omega_1 +s\delta = (\ell-m)\Lambda_0+m\Lambda_1+s\delta$ for some  $0\le m\le \ell, s\in\mathbb Z$. Note that, if $\ell=0$, then $\Lambda=s\delta$ and $\widehat{\mathcal W}\Lambda = \{\Lambda\}$. Thus, we assume from now on that $\ell>0$. One then easily sees that
\begin{gather}\notag
\sigma_k\Lambda = \ell\Lambda_0-(2k\ell+m)\omega_1-(k(k\ell+m)-s)\delta\\ \label{e:sl2orbit} \quad\text{and}\quad \\ \notag
\sigma_{k}s_1\Lambda = \ell\Lambda_0-(2k\ell-m)\omega_1-(k(k\ell-m)-s)\delta
\end{gather}
for all $k\in\mathbb Z$. In particular,  ${\widehat{\mathcal W}_\Lambda}^- = \{\sigma_k,\sigma_{k'}s_1:k\ge 0,k'\ge 1-\delta_{m,0}\}$.

\subsection{The Proof of Proposition \ref{c:sl2l}}\label{ss:sl2l}
We write down the proof for $i=0$. The case $i=1$ easily follows from that one together with \eqref{e:multrelnot0} and, thus, we omit the details. Setting $V=V(\Lambda_0)\otimes V(\Lambda)$ and $m=m_0=\Lambda(h_1)$,
we want to show that
\begin{equation}\label{e:sl2l0}
[V:V(\Phi)] =  \sum_{ (\lambda,r)\in\Gamma_{\Phi}}\ \lim_{k\to\infty}\  \alpha_\lambda^\ell \left(k\ell+\frac{m-\lambda}{2}\ ,\  r+ k(k\ell+m)-s\right).
\end{equation}	
	To shorten notation, set $\ell'=\ell+1$. It follows from \eqref{e:sl2orbit} that 
	\begin{equation*}
	V_{\sigma_k}(\Lambda) \cong D(\ell,2k\ell+m, -k(k\ell+m)+s) \quad\text{for}\quad k\ge 0.
	\end{equation*}
	Observe that 
	\begin{equation*}
	\sigma_k\preccurlyeq \sigma_{k+1}s_1\preccurlyeq\sigma_{k+1} \quad\text{for all}\quad k\ge 0.
	\end{equation*}
	It then follows from Theorem \ref{t:multrel} together with \eqref{e:BruhatDem}  that
	\begin{align*}
	[V:V(\Phi)] & =  \sum_{ (\lambda,r)\in\Gamma_\Phi}\ \lim_{k\to\infty}\ [D(\ell,2k\ell+m, -k(k\ell+m)+s): D(\ell',\lambda,r)]\\
	& =  \sum_{ (\lambda,r)\in\Gamma_\Phi}\  \lim_{k\to\infty}\ [D(\ell,2k\ell+m): D(\ell',\lambda,r+k(k\ell+m)-s)]\\
	& \stackrel{\eqref{e:defalfa}}{=}  \sum_{ (\lambda,r)\in\Gamma_\Phi}\ \lim_{k\to\infty}\ \alpha_\lambda^\ell \left(\frac{2k\ell+m-\lambda}{2},r+k(k\ell+m)-s\right),
	\end{align*}
	thus proving \eqref{e:sl2l0}. 

\subsection{The Proof of Proposition \ref{t:sl21}}\label{ss:sl21} 
Recall that $V=V(\Lambda_0)\otimes V(\Lambda_i)$. It follows from Proposition \ref{c:sl2l} with $\Lambda=\Lambda_i$ that
\begin{align}\label{e:sl2l=1}
[V:V(\Phi)] & =  \sum_{ (\lambda,r)\in\Gamma_\Phi}\ \lim_{k\to\infty}\  
\alpha_\lambda^1 \left(k -\frac{\lambda-i}{2},r+k(k+i)\right).
\end{align}	
Fix $\Phi$ such that $[V:V(\Phi)]\ne 0$. Evidently, $\Phi\le\Lambda_0+\Lambda_i$ and, hence, $\Phi = 2\Lambda_0+m\omega_1-s\delta$ for some $m\in\{0,1,2\}$ and $s\in\mathbb Z_{\ge 0}$. 
By \eqref{e:sl2orbit}, 
\begin{equation}\label{e:sl2orbitl=2}
\Gamma_\Phi = \{(4l+m, -l(2l+m)-s):l\ge 0\}\ \cup\ \{(4l-m, -l(2l-m)-s):l\ge 1\}.
\end{equation}
In particular,
\begin{equation*}
(\lambda,r)\in\Gamma_\Phi \quad\Rightarrow\quad \lambda-m\in 2\mathbb Z.
\end{equation*}
Since,  $[V:V(\Phi)]\ne 0$, \eqref{e:sl2l=1} implies that there must exist $(\lambda,r)\in\Gamma_\Phi$ such that $\lambda-i\in 2\mathbb Z$ and, therefore, $m-i\in 2\mathbb Z.$ 
This proves that $m=2j+i$ and the condition $0\le m\le 2$ implies that $0\le j\le \delta_{i,0}$ as claimed. Moreover, it also follows that
\begin{equation*}
(\lambda,r)\in\Gamma_\Phi \quad\Rightarrow\quad \lambda-i\in 2\mathbb Z_{\ge0}.
\end{equation*}

It was shown in \cite{cssw:demflag} that
\begin{equation}\label{e:csswmf2}
[W(\mu):D(2,\lambda)](q) = 
\begin{cases}
q^{p\lceil \mu/2\rceil} \tqbinom{\lfloor \mu/2\rfloor}{p}_q, & \text{if } p:=\frac{\mu-\lambda}{2}\in \mathbb Z_{\ge 0},\\
0, &\text{otherwise.}
\end{cases}
\end{equation}
Here,
\begin{equation*}
\qbinom{m}{p}_q = \prod_{n=0}^{p-1} \frac{1-q^{m-n}}{1-q^{p-n}} \quad\text{for} \quad m,p\in\mathbb Z_{\ge 0}, \ p\le m.
\end{equation*}
Given $m,p\in\mathbb Z_{\ge 0}, 0\le p\le m$, define $\beta^\pm_{m,m-p}(r)\in\mathbb Z$ by 
\begin{equation*}
q^{(m+\delta_{1,\pm 1})p}\qbinom{m}{p}_q = \sum_{r\in\mathbb Z} \beta^\pm_{m,m-p}(r)\ q^r.
\end{equation*}
Set $\beta^\pm_{m,l}(r)=0$ if $l\ne m-p$ for some $p$ as above. 
In particular,   \eqref{e:csswmf2} can be rewritten as 
\begin{equation}\label{e:csswmf2m}
[W(\mu):D(2,\lambda,r)] = \begin{cases}
\beta^{{\rm sign}(-1)^{\mu+1}}_{\lfloor \mu/2\rfloor,\lfloor \lambda/2\rfloor}(r), & \text{ if } \mu-\lambda\in2\mathbb Z_{\ge 0},\\
0,& \text{ otherwise,}
\end{cases}
\end{equation}
and
\begin{equation*}
\alpha_\lambda^1(m,r) = [W(\lambda+2m):D(2,\lambda,r)] = \beta^{{\rm sign}(-1)^{\lambda+1}}_{\lfloor \lambda/2\rfloor+m,\lfloor \lambda/2\rfloor}(r) 
\end{equation*}
for all $\lambda,m,r\in\mathbb Z, \lambda\ge 0$.
Plugging back in \eqref{e:sl2l=1}, we get
\begin{align}\label{e:sl2beta}
[V:V(\Phi)] = \sum_{ (\lambda,r)\in\Gamma_\Phi}\ \lim_{k\to\infty}\  
\beta^{{\rm sign}(-1)^{i+1}}_{k,\lfloor \lambda/2\rfloor}(r+k(k+i)).
\end{align}	
	
Before continuing, we recall that the coefficient of $q^r$ in $\tqbinom{m}{p}_q$ is equal to $\bro_{m-p}^p(r)$,
the number of partitions of $r$ in at most $p$ parts all bounded by $m-p$ \cite[Chapter 3]{andrews}. Therefore, 
\begin{equation}\label{e:betaro}
\beta^-_{k,k-p}(r)= \bro_{k-p}^p(r-kp) \qquad\text{and}\qquad \beta^+_{k,k-p}(r)= \bro_{k-p}^p(r-(k+1)p)
\end{equation}
for all $0\le p\le k, r\in\mathbb Z$.

Assume first that $i=0$. Then, $j\in\{0,1\}, \Phi=2\Lambda_j-s\delta$, and  \eqref{e:sl2orbitl=2} is equivalent to:
	\begin{equation*}
	(\lambda,r)\in\Gamma_\Phi \quad\Leftrightarrow\quad \lambda=2(2l+j) \quad\text{and}\quad r= -2l(l+j)-s \quad\text{for some}\quad l\ge 0.
	\end{equation*}
	Plugging this back in \eqref{e:sl2beta}, we get
	\begin{align*}
	[V:V(\Phi)] & =  \sum_{l\ge 0}\ \lim_{k\to\infty}\ 	\beta^-_{k,2l+j}\left(k^2-2l(l+j)-s\right).
	\end{align*}	
Thus, we are left to show that
\begin{gather}\notag
\lim_{k\to\infty}\ 	\beta^-_{k,2l+j}\left(k^2-2l(l+j)-s\right) =  \bro_{2l+j}(s-2l^2-2jl-j)\\ \text{and}\\ \notag 
s-2l^2-2jl-j < 0 \qquad\text{if}\qquad s<j \quad\text{or}\quad l> \frac{-j+\sqrt{2s-j}}{2}.
\end{gather}
The second of these statements follows from a simple analysis of the quadratic polynomial in $l$ 
$$f_{j,s}(l)= s-j^2-2jl-2l^2,$$
observing that $j^2=j$.  
To prove the first statement, given $l\ge 0$, let $k\ge 2l+j$ and note that \eqref{e:betaro} implies
\begin{equation}
\beta^-_{k,2l+j}(k^2-2l(l+j)-s) = \bro_{b}^{k-b}(a_k) \quad\text{where}\quad b = 2l+j,\ a_k=kb-2l(l+j)-s.
\end{equation}
Observe that
\begin{equation*}
a_k = b(k-b)- f_{j,s}(l)
\end{equation*}
and, hence,
\begin{equation*}
f_{j,s}(l)<0 \quad\Rightarrow\quad \bro_{b}^{k-b}(a_k)=0.
\end{equation*}
If $f_{j,s}(l)\ge 0$ and $k$ is sufficiently large, i.e., if $k\ge f_{j,s}(l)+b$, all partitions of $f_{j,s}(l)$ have at most $k-b$ parts and, hence, we have a bijection
\begin{equation*}
\mathscr P_b(f_{j,s}(l)) \to \mathscr P_b^{k-b}(a_k), \qquad (\lambda_1\ge\lambda_2\ge \dots\ge \lambda_{k-b}\ge 0) \mapsto (b-\lambda_{k-b}\ge \dots\ge b-\lambda_1\ge 0),
\end{equation*}
where $\mathscr P_b(a)$ is the set of partitions of $a$ whose parts are bounded by $b$ and $\mathscr P_b^p(a)$ is the subset of partitions with at most $p$ parts. Thus,
\begin{equation*}
\beta^-_{k,2l+j}(k^2-2l(l+j)-s) = \bro_{b}^{k-b}(a_k) = \bro_b(f_{j,s}(l)) \quad\text{for all}\quad k\gg 0,
\end{equation*}
proving Theorem \ref{t:sl21} for $i=0$.
	
If $i=1$, then $\Phi=\Lambda_0+\Lambda_1-s\delta$ and \eqref{e:sl2orbitl=2} becomes
\begin{equation*}
\Gamma_\Phi  = \left\{(4l+ 1,-l(2l+ 1)-s):l\ge 0\right\}\ \cup\ \left\{(4l- 1,-l(2l- 1)-s): l\ge 1\right\}.
\end{equation*}
Thus, \eqref{e:sl2beta} becomes
\begin{align}\label{e:sl2beta1}
[V:V(\Phi)] & =  \sum_{\epsilon=\pm 1}\ \sum_{l\ge \delta_{\epsilon,-1}}\ \lim_{k\to\infty}\ \beta^+_{k,2l-\delta_{\epsilon,-1}}\left(k(k+\epsilon)-l(2l+\epsilon)-s\right).
\end{align}
Using \eqref{e:betaro} we see that
\begin{equation*}
\beta^+_{k,2l-\delta_{\epsilon,-1}}\left(k(k+\epsilon)-l(2l+\epsilon)-s\right) = \bro_{b_\epsilon}^{k-b_\epsilon}(a_{\epsilon,k})
\end{equation*}
where  
\begin{equation*}
b_\epsilon = 2l-\delta_{\epsilon,-1},\quad a_{\epsilon,k}=b_\epsilon(k-b_\epsilon)- f_{\epsilon,s}(l), \quad\text{and}\quad f_{\epsilon,s}(l)=-2l^2-\epsilon l +s.
\end{equation*}
Recall the following identity \cite[Equation (3.2.6)]{andrews} 
\begin{equation*}
\bro_{b-1}^p(r-p) = \bro_b^p(r)-\bro_b^{p-1}(r)
\end{equation*}
and observe that, for fixed $l\ge 1$,
\begin{equation*}
b_{-1} = b_1-1 \qquad\text{and}\qquad a_{-1,k} = a_{1,k} -(k-2l+1).
\end{equation*}
It follows that
\begin{equation*}
\bro_{b_{-1}}^{k-b_{-1}}(a_{-1,k}) = \bro_{b_{1}-1}^{k-b_{-1}}(a_{1,k} - (k-b_{-1})) = \bro_{b_{1}}^{k-b_{-1}}(a_{1,k})  - \bro_{b_{1}}^{k-b_{1}}(a_{1,k}).
\end{equation*}
Hence, 
\begin{equation*}
\beta^+_{k,2l}\left(k(k+1)-l(2l+1)-s\right) + \beta^+_{k,2l-1}\left(k(k-1)-l(2l-1)-s\right) = \bro_{b_{1}}^{k-b_{-1}}(a_{1,k})
\end{equation*}
for all $l\ge 1$. Plugging this back in \eqref{e:sl2beta1} we get 	
\begin{align*}
[V:V(\Phi)] & = \lim_{k\to\infty}\ \bro^{k}_{0}(-s) \ +\ \sum_{l\ge 1}\ \lim_{k\to\infty}\ \bro_{b_{1}}^{k-b_{-1}}(a_{1,k}) = \delta_{s,0}\ +\ \sum_{l\ge 1}\ \lim_{k\to\infty}\ \bro_{b_{1}}^{k-b_{-1}}(a_{1,k}).
\end{align*}
Since
\begin{equation*}
a_{1,k} = b_1(k-b_{-1}) - f_{-1,s}(l) \quad\text{for}\quad l\ge 1,
\end{equation*}
proceeding as before, we see that 
\begin{equation*}
\bro_{b_{1}}^{k-b_{-1}}(a_{1,k}) = \bro_{b_1}(f_{-1,s}(l)) \quad\text{for all}\quad k\gg 0 \qquad\text{and}\qquad f_{-1,s}(l) < 0 \quad\text{if}\quad l> L.
\end{equation*}
Since $\bro_{0}(f_{-1,s}(0))=\delta_{s,0}$, this completes the proof of the proposition.

\bibliographystyle{amsplain}

\end{document}